\documentclass[11pt]{amsart}
\usepackage{latexsym}
\usepackage{amsmath}
\usepackage{amssymb}

\newcommand{\eproof}{\mbox{\ }\hfill $\Box$ \par \vskip 10pt}

\newtheorem{Theorem}{Theorem}[section]
\newtheorem{lemma}[Theorem]{Lemma}
\newtheorem{prop}[Theorem]{Proposition}

\numberwithin{equation}{section}

\def\cal{\mathcal}

\begin{document}

\title[Transmission eigenvalue-free regions near the real axis]{Transmission eigenvalue-free regions near the real axis
in the anisotropic case}

\author[G. Vodev]{Georgi Vodev}

\address {Universit\'e de Nantes, Laboratoire de Math\'ematiques Jean Leray, 2 rue de la Houssini\`ere, BP 92208, 44322 Nantes Cedex 03, France}
\email{Georgi.Vodev@univ-nantes.fr}

\date{}

\begin{abstract} 
We consider the anisotropic 
interior transmission problem with one complex-valued refraction index, that is, with a damping term
which does not vanish on the boundary. 
Under the condition that all geodesics reach the boundary, 
for a class of strictly concave domains, we obtain a transmission eigenvalue-free region of the form 
 $\{\lambda\in \mathbb{C}:C_N(|\lambda|+1)^{-N}\le|{\rm Im}\,\lambda|\le C(|\lambda|+1)^{-1}\}$, where 
 $N>1$ is arbitrary. Under extra conditions on the coefficients we get a larger 
 transmission eigenvalue-free region of the form 
 $\{\lambda\in \mathbb{C}:C_N(|\lambda|+1)^{-N}\le|{\rm Im}\,\lambda|\le C\}$.

\end{abstract} 

\maketitle

\setcounter{section}{0}
\section{Introduction}

Let $\Omega\subset\mathbb{R}^d$, $d\ge 2$, be a bounded, connected domain with a $C^\infty$ smooth boundary $\Gamma=\partial\Omega$.
Let $f\in C^\infty(\mathbb{R}^d)$ be such that $f<0$ in $\Omega$, $f>0$ in $\mathbb{R}^d\setminus\Omega$, $f=0$, $df\neq 0$ on $\Gamma$.
Let $g\in C^\infty(T^*\Omega)$ be a Hamiltonian of the form
$$g(x,\xi)=\sum_{i,j=1}^dg_{ij}(x)\xi_i\xi_j\ge C|\xi|^2,\quad C>0.$$
The boundary $\Gamma$ will be said to be $g-$strictly concave (viewed from the interior) iff for any $(x,\xi)$
satisfying 
$$f(x)=0, g(x,\xi)=1, \{g,f\}(x,\xi)=0,$$
we have 
$$\{g,\{g,f\}\}(x,\xi)>0,$$
where $\{\cdot,\cdot\}$ denotes the Poisson brackets. We also define the bicharacteristic flow 
$$\Phi(t):T^*\Omega\to T^*\Omega,\quad t\in \mathbb{R},$$
associated to $g$ by $\Phi(t)(x^0,\xi^0)=(x(t),\xi(t))$, where $(x(t),\xi(t))$ solves the Hamilton equation
\begin{equation}\label{eq:1.1}
\frac{\partial x(t)}{\partial t}=\frac{\partial g(x,\xi)}{\partial\xi},\quad 
\frac{\partial\xi(t)}{\partial t}=-\frac{\partial g(x,\xi)}{\partial x},\quad x(0)=x^0,\,\xi(0)=\xi^0.
\end{equation}
The solutions to (\ref{eq:1.1}) exist for $t\in I:=(t^-,t^+)$, where $0<\pm t^\pm\le \infty$ depends on $(x^0,\xi^0)$
and is such that either 
$\pm t^\pm<\infty$ and 
$$x(t^\pm):=\lim_{t\to t^\pm,\,t\in I}x(t)\in\Gamma,$$ 
or $t^\pm=\pm\infty$. Clearly, if $(x^0,\xi^0)\in{\cal G}:=\{(x,\xi)\in T^*\Omega: g(x,\xi)=1\}$, then
$(x(t),\xi(t))\in{\cal G}$ for all $t$ for which the solution exists. 
The curves $(x(t),\xi(t))$ in ${\cal G}$ are called null 
bicharacteristics of the Hamiltonian $g(x,\xi)-1$. The projection in the $x-$space of a null bicharacteristic
is called a $g-$geodesic, that is, the curve $x(t)$ in $\Omega$. 
 
Consider the interior transmission problem

\begin{equation}\label{eq:1.2}
\left\{
\begin{array}{l}
(\nabla c_1(x)\nabla+\lambda^2n_1(x)+i\lambda m(x))u_1=0\quad \mbox{in}\quad\Omega,\\
(\nabla c_2(x)\nabla+\lambda^2n_2(x))u_2=0\quad \mbox{in}\quad\Omega,\\
u_1=u_2,\,c_1\partial_\nu u_1=c_2\partial_\nu u_2 \quad\mbox{on}\quad\Gamma,
\end{array}
\right.
\end{equation}
where $\lambda\in \mathbb{C}$, $\nu$ denotes the Euclidean unit inner normal to $\Gamma$ and $c_j,n_j,m\in C^\infty(\overline\Omega)$, $j=1,2$, are real-valued functions satisfying $c_j(x)>0$, $n_j(x)>0$. We also suppose that either $m(x)\ge 0$ for all $x\in \overline\Omega$
or $m(x)\le 0$ for all $x\in \overline\Omega$. 
If the equation (\ref{eq:1.2})
has a non-trivial solution $(u_1,u_2)$ the complex number $\lambda$ is said to be an interior transmission eigenvalue. 
We refer to the very nice survey \cite{kn:CCH} where the role the transmission eigenvalues play in the scattering theory and the inverse
problems is explained and many properties are described. 

We suppose that the coefficients satisfy the condition
\begin{equation}\label{eq:1.3}
c_1(x)\neq c_2(x) \quad\mbox{on}\quad\Gamma.
\end{equation}
Then the interior transmission problem (\ref{eq:1.2}) is called anisotropic. 
When $m\equiv 0$ it is known (e.g. see \cite{kn:NN1}, \cite{kn:NN2}, \cite{kn:LV1}) that the transmission eigenvalues form a discreet set in $\mathbb{C}$
with no finite accumulation points under the conditions (\ref{eq:1.3}) and
\begin{equation}\label{eq:1.4}
c_1(x)n_1(x)\neq c_2(x)n_2(x) \quad\mbox{on}\quad\Gamma.
\end{equation}
This is also true in the isotropic case ($c_1\equiv c_2\equiv 1$) again under the condition (\ref{eq:1.4})
(e.g. see \cite{kn:NN1}, \cite{kn:FN}, \cite{kn:LV2}, \cite{kn:Sy}). Note that the condition (\ref{eq:1.4}) 
is crucial to prove discreetness of the transmission eigenvalues 
as well as to study their location on the complex plane 
in both isotropic and anisotropic cases when $m$ is identically zero. When $m$ is not identically zero this seems not to be the case. 

  It is well-known that the counting function of the transmission eigenvalues has Weyl type
asymptotics. Such asymptotics are proven in \cite{kn:PV}, \cite{kn:R2} when the coefficients are $C^\infty$
smooth and in \cite{kn:NN2}, \cite{kn:FN} when the coefficients have very low regularity. Note that the result in 
\cite{kn:R2} also holds when the function $m$ is not identically zero. The proof in \cite{kn:PV}
relies heavily on the fact that there are no transmission eigenvalues outside some parabolas.
In the case of $C^\infty$ coefficients, parabolic transmission eigenvalue-free regions were obtained in
\cite{kn:V1} and latter on improved in \cite{kn:V2}. The main ingridient in the proof in
\cite{kn:V1} is the construction of a semiclassical parametrix for the Dirichlet-to-Neumann map valid in the region
$$\{\lambda\in \mathbb{C}:|{\rm Im}\,\lambda|\ge (|\lambda|+1)^{1/2+\epsilon},\,|\lambda|\gg 1\},$$
 $0<\epsilon\ll 1$ being arbitrary.
Then the problem of getting transmission eigenvalue-free regions is reduced to that one of inverting 
semiclassical pseudo-differential operators. The parametrix construction in
\cite{kn:V1} requires solving the eikonal and the transport equations modulo some error terms, which is particularly
hard to do when the boundary data is microlocally supported near the glancing region. In contrast, the
parametrix construction is relatively
easy to do when the boundary data is microlocally supported in the elliptic region. In this case the 
parametrix construction is valid in the region
$$\{\lambda\in \mathbb{C}:|{\rm Im}\,\lambda|\ge (|\lambda|+1)^{-N},\,|\lambda|\gg 1\},$$
 $N\gg 1$ being arbitrary (see Section 5).
When the boundary data is microlocally supported in the hyperbolic region, one can construct a parametrix 
valid in the region
$$\{\lambda\in \mathbb{C}:|{\rm Im}\,\lambda|\ge C,\,|\lambda|\gg 1\},$$
 $C\gg 1$ being a suitably choosen constant (see Section 4
of \cite{kn:V2}).
Most probably the results in \cite{kn:V1} still hold 
when the function $m$ is not identically zero, since the term $i\lambda m$ is a lower-order perturbation
which does not affect the parametrix construction significantly. In particular, the eikonal equation does not depend on $m$.
Note also that parabolic transmission eigenvalue-free regions have been recently obtained in
\cite{kn:V3} for coefficients of low regularity.

When the function $m$ is not identically zero, it is well-known that there are no real transmission eigenvalues. 
Moreover, using the interpolation estimates
in \cite{kn:LR}, Robbiano showed (see Theorem 8 of \cite{kn:R1}) that there are no 
transmission eigenvalues in the region
\begin{equation}\label{eq:1.5}
|{\rm Im}\,\lambda|\le C_1e^{-C_2|\lambda|},\quad C_1,C_2>0.
\end{equation}
In fact, Robbiano considered the isotropic interior transmission problem, that is, when $c_1\equiv c_2\equiv 1$, 
but his arguments extend to the anisotropic one as well.

Our goal in the present paper is to study the location of the transmission eigenvalues near the real axis under the stronger condition 
\begin{equation}\label{eq:1.6}
m(x)\neq 0\quad\mbox{on}\quad\Gamma.
\end{equation}
We impose some {\it nontrapping} conditions, namelly that all geodesics associated to the Hamiltonians
$g_j=\frac{c_j(x)}{n_j(x)}|\xi|^2$, $j=1,2$, reach the boundary in a finite time.
More precisely, we make the following assumption:
\begin{equation}\label{eq:1.7}
\mbox{There exists a constant $T>0$ such that for any $g_j-$geodesic $\gamma(t)$ with $\gamma(0)\in\Omega$}$$
$$\mbox{ there are $t_\gamma^\pm\in \mathbb{R}$,
$0\le\pm t_\gamma^\pm\le T,$ such that $\gamma(t_\gamma^\pm)\in \Gamma$.}
\end{equation}
Note that the condition (\ref{eq:1.7}) is trivially fulfilled if the functions $\frac{c_j(x)}{n_j(x)}$, $j=1,2$, 
are constants, since in this case
the $g_j-$geodesics are lines. 
 Our main result is the following

\begin{Theorem} Suppose that the conditions (\ref{eq:1.3}), (\ref{eq:1.6}) and (\ref{eq:1.7}) are satisfied. Suppose also that 
$\Gamma$ is $g_2-$strictly concave. Then for every $N>1$ there exist constants $C, C_N>0$,
$C$ being independent of $N$, such that there are no transmission eigenvalues in the region
  \begin{equation}\label{eq:1.8}
C_N(|\lambda|+1)^{-N}\le|{\rm Im}\,\lambda|\le C(|\lambda|+1)^{-1}.
\end{equation}
If in addition we assume the condition 
 \begin{equation}\label{eq:1.9}
\frac{n_1(x)}{c_1(x)}>\frac{n_2(x)}{c_2(x)},\quad c_1(x)<c_2(x),\quad\mbox{on}\quad\Gamma,
\end{equation}
then there are no transmission eigenvalues in the region
  \begin{equation}\label{eq:1.10}
C_N(|\lambda|+1)^{-N}\le|{\rm Im}\,\lambda|\le C.
\end{equation}
 \end{Theorem}

\noindent
{\bf Remark 1.} The fact that we can take an arbitrary $N$ comes from the $C^\infty-$ smoothness of the boundary $\Gamma$ and
the coefficients $c_j$, $n_j$, $m$ near $\Gamma$. Therefore, one can expect that if more regularity is assumed (e.g. Gevrey class
or analyticity), larger eigenvalue-free regions exist. Indeed, using the pseudo-differential calculus with analytic symbols
developed in \cite{kn:Sj}, one can replace $C_N(|\lambda|+1)^{-N}$ in (\ref{eq:1.8}) and (\ref{eq:1.10}) by
$\beta_1e^{-\beta_2|\lambda|}$, $\beta_1, \beta_2>0$ being constants, provided $\Gamma$ and
the coefficients $c_j$, $n_j$, $m$ near $\Gamma$ are supposed analytic.

\noindent
{\bf Remark 2.} We expect that the region 
$$|{\rm Im}\,\lambda|\le C_N(|\lambda|+1)^{-N}$$
is also free of transmission eigenvalues, but the proof would be harder since for such $\lambda$ the parametrix construction 
(in the case $m\equiv 0$) of 
the interior Dirichlet-to-Neumann map in the elliptic region is not valid (see Section 5). Therefore, to prove such a statement a different
approach is required.

\noindent
{\bf Remark 3.} We do not need the condition (\ref{eq:1.4}) in the proof of the above theorem because of  
the condition (\ref{eq:1.6}) which is strong enough to replace it.

\noindent
{\bf Remark 4.} We expect that smaller transmission eigenvalue-free regions still exist without the strict concavity assumption.
To get such regions one needs to prove an analog of Theorem 3.1 under the nontrapping condition, only.

\noindent
{\bf Remark 5.}
It is worth noticing that free regions similar to (\ref{eq:1.10}) exist for the resonances associated to exterior transmission problems under conditions
on the coefficients similar to (\ref{eq:1.9}) (see \cite{kn:CPV}, \cite{kn:G}). The big difference with the problem studied in the present paper
is that, when the domain is strictly convex, the exterior Dirichlet-to-Neumann map has a complete parametrix valid in any strip
$|{\rm Im}\,\lambda|\le M$, $|\lambda|\ge C_M\gg 1$, $M>0$ being arbitrary. As mentioned above, this is no longer true for 
the interior Dirichlet-to-Neumann map, which makes the study of the location of the transmission eigenvalues in a strip near the real axis
much harder than the study of the location of the resonances. 

To prove Theorem 1.1 we need to show that if $\lambda$ belongs to the eigenvalue-free regions, then the solution $(u_1,u_2)$ to the equation 
(\ref{eq:1.2}) is identically zero. Clearly, it suffices to show that the function $f=u_1|_\Gamma=u_2|_\Gamma$ is identically zero.
To this end, we use that $T(\lambda)f\equiv 0$, where 
$$T(\lambda)=c_1{\cal N}_1(\lambda)-c_2{\cal N}_2(\lambda),$$ 
${\cal N}_1(\lambda)$ and ${\cal N}_2(\lambda)$ being the corresponding interior Dirichlet-to-Neumann maps
(see Section 6). Roughly speaking, the idea
is to bound from above the $H^1(\Gamma)$ norm of $f$ by the same norm times a factor of the form $\tau(|\lambda|)|{\rm Im}\,\lambda|^{1/2}$. 
Then we would have $f\equiv 0$ in the region $|{\rm Im}\,\lambda|\le (2\tau(|\lambda|))^{-2}$. As shown in Section 5, the operators ${\cal N}_1(\lambda)$ and ${\cal N}_2(\lambda)$ can be approximated in the corresponding elliptic regions ${\cal E}_1$ and ${\cal E}_2$
by $h-\Psi$DOs, $h=|\lambda|^{-1}$, with symbols belonging to $S^1(\Gamma)$, 
provided 
$$|{\rm Im}\,\lambda|\ge |\lambda|^{-N}, \quad |\lambda|\gg 1.$$  
Thus we conclude that the operator $T(\lambda)$ can be approximated in ${\cal E}_1\cap{\cal E}_2$
by a $h-\Psi$DO, say ${\cal T}$, with principal symbol 
$$\sigma_p({\cal T})=(c_2-c_1)r_0^{1/2}+{\cal O}(r_0^{-1/2}),\quad r_0\gg 1.$$
Hereafter $r_0$ denotes the principal symbol
of the positive Laplace-Beltrami operator on $\Gamma$. In other words, the condition (\ref{eq:1.3}) guarantees that
${\cal T}$ is an elliptic first order $h-\Psi$DO in the deep elliptic region $\{r_0\gg 1\}$. 
However, it does not guarantee that 
${\cal T}$ is elliptic everywhere in ${\cal E}_1\cap{\cal E}_2$. Nevertheless, the condition (\ref{eq:1.3}) implies that
the operator
${\cal T}:H^1(\Gamma)\to L^2(\Gamma)$
is invertible in the deep elliptic region. This fact allows us to conclude that the $H^1(\Gamma)$ norm of ${\rm Op}_h(1-\chi)f$
can be bounded by the $H^1(\Gamma)$ norm of $f$ times a small factor (see Lemma 6.2), where $\chi\in C^\infty(T^*\Gamma)$
is of compact support such that $1-\chi$ is supported in $\{r_0\gg 1\}$. On the other hand, the a priori estimates of the solution of
the Helmholtz equation with a damping term, obtained in Section 2 under the conditions (\ref{eq:1.6}) and (\ref{eq:1.7}),
allow us to control the $H^1(\Gamma)$ norm of $f$ by the $H^1(\Gamma)$ norm of ${\rm Op}_h(1-\chi)f$ (see the
estimate (\ref{eq:2.5})). Thus we get the eigenvalue-free region (\ref{eq:1.8}). The proof of (\ref{eq:1.10}) is similar
observing that the condition (\ref{eq:1.9}) implies ${\cal E}_1\subset{\cal E}_2$ and that ${\cal T}$ is elliptic 
in ${\cal E}_1$. This allows us to get a better control of the $H^1(\Gamma)$ norm of $f$ by the $H^1(\Gamma)$ norm of ${\rm Op}_h(1-\chi)f$, provided $\chi$ is properly chosen (see the estimate (\ref{eq:2.27}) and Lemma 6.3), and hence to get a larger eigenvalue-free region. 

The paper is organized as follows. In Section 2 we obtain some a priori estimates for the solutions of the interior Helmholtz equation
(see equation (\ref{eq:2.1}) below) under the conditions (\ref{eq:1.6}) and (\ref{eq:1.7}). In particular, (\ref{eq:1.7}) allows us
to control the global $L^2$ norm of the solution to (\ref{eq:2.1}) by the $L^2$ norm near the boundary (see Proposition 2.2,
which is also valid when $m\equiv 0$). This kind of estimates are not new and even more general result exist (e.g. see Theorem 3.3
of \cite{kn:BLR}), but here we give a short proof which allows to see the role of the function $m$ and why the fact that $m$ does not change sign
matters. In Section 3 we obtain a priori estimates for the solutions of the interior Helmholtz equation with $m\equiv 0$ 
(see equation (\ref{eq:3.1}) below) under the condition (\ref{eq:1.7}) and that the boundary is $g_2-$strictly concave.
These conditions allow a good control of the global $L^2$ norm of the solution to (\ref{eq:3.1}) by the $L^2$ norm of the
Diriclet and Neumann data on the boundary (see Theorem 3.1). The strict concavity allows to get nice a priori estimates
near the boundary (see Proposition 3.2) which combinned with Proposition 2.2 imply Theorem 3.1. Note that 
Proposition 3.2 is proved in \cite{kn:CPV} and here we omit the proof since it is quite long. In Section 4 we define
the interior Dirichlet-to-Neumann map and prove some important properties. In Section 5 we recall the parametrix construction
of the interior Dirichlet-to-Neumann map in the elliptic region, which is well-known when $m\equiv 0$ (see \cite{kn:V1}, \cite{kn:V2}).
As mentioned above, the presence of the function $m$ does not change anything in the construction. The only difference is that 
when $m$ is not identically zero the parametrix construction is valid in the region 
$\{|{\rm Im}\,\lambda|\le C,\,|\lambda|\gg 1\}$ with a suitably choosen constant
$C>0$, while in the case when $m$ is identically zero it is valid in the region 
$\{|\lambda|^{-N}\le|{\rm Im}\,\lambda|\le M,\,|\lambda|\gg 1\}$, $N,M\gg 1$ being arbitrary. 
Finally, in Section 6 we prove
Theorem 1.1 by combining the a priori estimates from the previous sections. 

\section{A priori estimates under the condition (\ref{eq:1.6})}

In this section we will obtain a priori estimates for the solution to the equation
\begin{equation}\label{eq:2.1}
\left\{
\begin{array}{l}
(\nabla c(x)\nabla+\lambda^2n(x)+i\lambda m(x))u=\lambda v\quad \mbox{in}\quad\Omega,\\
u=f\quad\mbox{on}\quad\Gamma,\\
\end{array}
\right.
\end{equation}
where $\lambda\in\Lambda_k:=\{\lambda\in\mathbb{C}, |{\rm Im}\,\lambda|\le k,\,{\rm Re}\,\lambda\ge 1\}$, $k>0$, and  
$c,n,m\in C^\infty(\overline\Omega)$, $c>0$, $n>0$. The function $m$ is as in Section 1. 
Set $h=({\rm Re}\,\lambda)^{-1}$. 
 Throughout this paper $\|\cdot\|$, $\|\cdot\|_1$, $\|\cdot\|_0$ and $\|\cdot\|_{1,0}$ will denote the norms in $L^2(\Omega)$,
$H^1(\Omega)$, $L^2(\Gamma)$ and $H^1(\Gamma)$, respectively, where the norms in $H^1(\Omega)$ and $H^1(\Gamma)$ are the semiclassical ones, that is, 
$$\|u\|_1^2:=\sum_{0\le|\alpha|\le 1}\left\|(h\partial_x)^\alpha u\right\|^2,$$
$$\|u\|_{1,0}^2:=\sum_{0\le|\alpha|\le 1}\left\|(h\partial_x)^\alpha u\right\|_0^2.$$
In fact, throughout this paper all Sobolev spaces will be equipped with the semiclassical norm. 
Also, $\langle\cdot,\cdot\rangle$ and $\langle\cdot,\cdot\rangle_0$ will denote the scalar products in 
$L^2(\Omega)$ and $L^2(\Gamma)$, respectively. Given a symbol $a\in C^\infty(T^*\mathbb{R}^d)$, ${\rm Op}_h(a)$ will denote the
$h-\Psi$DO defined by
$$\left({\rm Op}_h(a)f\right)(x)=(2\pi h)^{-d}\int_{T^*\mathbb{R}^d} e^{i\langle x-y,\xi\rangle/h}
a(x,\xi)f(y)d\xi dy.$$
Similarly, given a symbol $\chi\in C^\infty(T^*\Gamma)$, ${\rm Op}_h(\chi)$ will denote the
$h-\Psi$DO defined by
$$\left({\rm Op}_h(\chi)f\right)(x)=(2\pi h)^{1-d}\int_{T^*\Gamma}e^{i\langle x-y,\xi\rangle/h}
\chi(x,\xi)f(y)d\xi dy.$$
Given any $\ell\in \mathbb{R}$, $S^\ell$ will denote the set of all symbols $a$ satisfying
$$\left|\partial_x^\alpha\partial_{\xi}^\beta a\right|\le C_{\alpha,\beta}(|\xi|+1)^{\ell-|\beta|}$$
for all multi-indices $\alpha$ and $\beta$. 
Our goal in this section is to prove the following

\begin{Theorem} \label{2.1} Suppose that $m$ satisfies (\ref{eq:1.6}) and that the condition (\ref{eq:1.7})
is fulfilled with $g_j$ replaced by $g=\frac{c(x)}{n(x)}|\xi|^2$. Let $u\in H^2(\Omega)$ satisfy equation (\ref{eq:2.1}) and set 
 $\omega=h\partial_\nu u|_\Gamma$. Then there are constants $k,\lambda_0>0$ such that for all 
 $\lambda\in\Lambda_k$ with ${\rm Re\,\lambda}\ge\lambda_0$ we have the estimates
\begin{equation}\label{eq:2.2}
\|\omega\|_0+\|u\|_1\lesssim \|v\|+\|f\|_{1,0},
\end{equation}
\begin{equation}\label{eq:2.3}
\|u\|\lesssim \|v\|+|{\rm Im}\,\langle f,c\omega\rangle_0|^{1/2},
\end{equation}
\begin{equation}\label{eq:2.4}
\|u\|_1\lesssim \|v\|+h^{1/2}\|f\|_{1,0}+|{\rm Im}\,\langle f,c\omega\rangle_0|^{1/2}.
\end{equation}
Given any function $\chi\in C^\infty(T^*\Gamma)$ of compact support, independent of $\lambda$, we have the estimate
\begin{equation}\label{eq:2.5}
\|f\|_{1,0}\lesssim |\lambda|^{1/2}\|v\|+|\lambda|^{1/2}|{\rm Im}\,\langle f,c\omega\rangle_0|^{1/2}
+\|{\rm Op}_h(1-\chi)f\|_{1,0}.
\end{equation}
\end{Theorem}

{\it Proof.} It is easy to see that it suffices to prove the estimates (\ref{eq:2.2}), (\ref{eq:2.3}) and (\ref{eq:2.4}) for real $\lambda\ge\lambda_0$.
Indeed, if $u$ satisfies equation (\ref{eq:2.1}) with complex $\lambda\in\Lambda_k$, then $u$ satisfies equation (\ref{eq:2.1}) with
$\lambda$ replaced by ${\rm Re\,\lambda}$ and $v$ replaced by 
$$\widetilde v=\frac{\lambda}{{\rm Re\,\lambda}}v+\frac{{\rm Im\,\lambda}}{{\rm Re\,\lambda}}\left(({\rm Im\,\lambda}-2i{\rm Re\,\lambda})n+m\right)u.$$
Therefore, by (\ref{eq:2.2}) applied with $\lambda$ replaced by ${\rm Re\,\lambda}$ and $v$ replaced by $\widetilde v$ we get
\begin{equation}\label{eq:2.6}
\|\omega\|_0+\|u\|_1\lesssim \|v\|+\|f\|_{1,0}+(k+h)\|u\|,\quad \lambda\in\Lambda_k.
\end{equation}
Taking $k$ and $h$ small enough we can absorb the last term in the right-hand side of (\ref{eq:2.6}) and conclude that
(\ref{eq:2.2}) holds for $\lambda\in\Lambda_k$ as well. Clearly, the same argument applies to (\ref{eq:2.3}) and (\ref{eq:2.4}) as well. 
Thus, in what follows we will prove these estimates  
for real $\lambda\ge\lambda_0$. Then we have $h=\lambda^{-1}$. 

Given a parameter $0<\delta\ll 1$, independent of $\lambda$, set 
$\Omega_\delta=\{x\in\Omega:{\rm dist}(x,\Gamma)<\delta\}$.  
We will first prove the following

\begin{prop} \label{2.2} Suppose that all $g-$geodesics reach the boundary $\Gamma$ in a finite time. Then there is a constant 
$\lambda_0>0$ such that for $\lambda\ge\lambda_0$ we have the estimate
\begin{equation}\label{eq:2.7}
\|u\|\lesssim \|v\|+\|u\|_{L^2(\Omega_\delta)}.
\end{equation}
\end{prop}

{\it Proof.} Let $\phi\in C_0^\infty(\mathbb{R}^d)$ be independent of $\lambda$
and such that $\phi=1$ in 
$\Omega\setminus\Omega_{\delta/2}$ and supp$\,\phi\subset\Omega\setminus\Omega_{\delta/3}$. Set
$$P=h^2n(x)^{-1}\nabla c(x)\nabla+1+ih\widetilde m(x),$$
where $\widetilde m(x)=\frac{m(x)}{n(x)}$. Note that $g(x,\xi)-1$ is the principal symbol of $-P$.

Fix a point $(x^0,\xi^0)\in{\cal G}$, $x^0\in{\rm supp}\,\phi$, and choose a real-valued function $a(x,\xi)\in C_0^\infty(T^*\mathbb{R}^d)$, $0\le a\le 1$, such that
$a=1$ in a small neighbourhood of $(x^0,\xi^0)$ and $a=0$ outside another small neighbourhood of $(x^0,\xi^0)$. 
To simplify the notations, in what follows in the proof of Proposition 2.2, given any $s\in \mathbb{R}$, $s\neq 0$, we will denote by
$\|\cdot\|_s$ the semiclassical norm in the Sobolev space $H^s(\mathbb{R}^d)$. We will also denote by
$\|\cdot\|$ and $\langle\cdot,\cdot\rangle$ the norm and the scalar product in $L^2(\mathbb{R}^d;n(x)dx)$, respectively,
where $n$ is an extension to the whole $\mathbb{R}^d$ such that $n(x)\ge n_0$ for all $x$ with some constant
$n_0>0$.
Proposition 2.2 is a consequence of the following 

\begin{lemma} \label{2.3} Under the assumptions of Proposition 2.2, there is a constant
$\lambda_0>0$ such that for $\lambda\ge\lambda_0$ we have the estimate
\begin{equation}\label{eq:2.8}
\|{\rm Op}_h(a)(\phi u)\|\lesssim \lambda\|P(\phi u)\|_{-1}+h\|\phi u\|.
\end{equation}
\end{lemma}

{\it Proof.} Without loss of generality we may suppose that $m\ge 0$. Then for $t>0$ define the function
$a_t(x,\xi)\in C_0^\infty(T^*\mathbb{R}^d)$ by 
$$a_t(x,\xi)=a(\Phi(t)(x,\xi))=a(x(t),\xi(t)).$$
In view of the assumption, taking ${\rm supp}\,a$ small enough we can arrange that 
 $\pi_x{\rm supp}\,a_t\subset \Omega$
for $0\le t\le t_a$ and 
$\pi_x{\rm supp}\,a_{t_a}\cap{\rm supp}\,\phi_1=\emptyset$
for some $0<t_a\le T$, where $T>0$ is a constant indpendent of $a$, 
$\pi_x$ denotes the projection $(x,\xi)\to x$ and 
$\phi_1\in C_0^\infty(\mathbb{R}^d)$, independent of $\lambda$, is such that $\phi_1=1$ in 
$\Omega\setminus\Omega_{\delta/4}$ and supp$\,\phi_1\subset\Omega\setminus\Omega_{\delta/5}$. 
Denote by $P^*$ the formal adjoint operator of $P$ with respect to $\langle\cdot,\cdot\rangle$. We have 
$P-P^*=2ih\widetilde m.$

It is easy to check that the equation (\ref{eq:1.1}) implies the identity
$$\frac{\partial a_t}{\partial t}+\{g,a_t\}=0,\quad 0\le t\le t_a.$$
Hence the operator 
$$Q_t:=\lambda{\rm Op}_h(\partial_ta_t)+i\lambda^2[P,{\rm Op}_h(a_t)],\quad 0\le t\le t_a,$$
is a zero order $h-\Psi$DO, and hence uniformly bounded on $L^2(\mathbb{R}^d)$. We have
$$-\frac{1}{2}\frac{d}{dt}\|{\rm Op}_h(a_t)(\phi u)\|^2=-{\rm Re}\,\left\langle
{\rm Op}_h(\partial_ta_t)(\phi u),{\rm Op}_h(a_t)(\phi u)\right\rangle$$
$$=-\lambda{\rm Im}\,\left\langle
[P,{\rm Op}_h(a_t)](\phi u),{\rm Op}_h(a_t)(\phi u)\right\rangle$$
$$-h{\rm Re}\,\left\langle
Q_t(\phi u),{\rm Op}_h(a_t)(\phi u)\right\rangle$$
$$=\lambda{\rm Im}\,\left\langle
{\rm Op}_h(a_t)P(\phi u),{\rm Op}_h(a_t)(\phi u)\right\rangle$$ 
$$-2^{-1}\lambda{\rm Im}\,\left\langle
(P-P^*){\rm Op}_h(a_t)(\phi u),{\rm Op}_h(a_t)(\phi u)\right\rangle$$
$$-h{\rm Re}\,\left\langle
Q_t(\phi u),{\rm Op}_h(a_t)(\phi u)\right\rangle$$
$$=\lambda{\rm Im}\,\left\langle
{\rm Op}_h(a_t)P(\phi u),{\rm Op}_h(a_t)(\phi u)\right\rangle$$ 
$$-\left\langle 
\widetilde m{\rm Op}_h(a_t)(\phi u),{\rm Op}_h(a_t)(\phi u)\right\rangle$$
$$-h{\rm Re}\,\left\langle
Q_t(\phi u),{\rm Op}_h(a_t)(\phi u)\right\rangle$$
$$\lesssim  \left(\lambda\|{\rm Op}_h(a_t)P(\phi u)\|+
h\|\phi u\|\right)\|{\rm Op}_h(a_t)(\phi u)\|.$$
Since the operator ${\rm Op}_h(a_t):H^{-1}(\mathbb{R}^d)\to L^2(\mathbb{R}^d)$ is uniformly bounded, we have
$$-\frac{d}{dt}\|{\rm Op}_h(a_t)(\phi u)\|\lesssim\lambda\|P(\phi u)\|_{-1}+h\|\phi u\|.$$
Thus we get
$$\|{\rm Op}_h(a)(\phi u)\|-\|{\rm Op}_h(a_{t_a})(\phi u)\|=
-\int_0^{t_a}\frac{d}{dt}\|{\rm Op}_h(a_t)(\phi u)\|dt$$
$$\lesssim\lambda\|P(\phi u)\|_{-1}+h\|\phi u\|.$$
On the other hand, we have
$$\|{\rm Op}_h(a_{t_a})(\phi u)\|=\|[\phi_1,{\rm Op}_h(a_{t_a})](\phi u)\|\lesssim h\|\phi u\|,$$
where we have used that $\phi\equiv \phi_1\phi$ and $\phi_1a_{t_a}\equiv 0$. 
\eproof

Let $A\in C_0^\infty(T^*\mathbb{R}^d)$ be such that $A=1$ in a small neighbourhood of $\{(x,\xi)\in{\cal G}: x\in {\rm supp}\,\phi_1\}$
and $A=0$ outside another small neighbourhood. Then the operator $P$ is elliptic on supp$(1-A)$, that is, the principal symbol of $P$
satisfies $|\sigma_p(P)|\ge C\langle\xi\rangle^2$ on supp$(1-A)$, where $C>0$ is some constant. Hence
\begin{equation}\label{eq:2.9}
\|{\rm Op}_h(1-A)(\phi u)\|\lesssim \|P(\phi u)\|_{-1}+h\|\phi u\|$$
$$\lesssim \|\phi Pu\|+\|[P,\phi]\widetilde \phi u\|_{-1}+h\|\phi u\|$$
$$\lesssim \|\phi Pu\|+\|\widetilde \phi u\|+h\|\phi u\|$$
$$\lesssim \|Pu\|+\|u\|_{L^2(\Omega_\delta)}+h\|\phi u\|,
\end{equation}
where $\widetilde\phi\in C_0^\infty(\mathbb{R}^d)$ is such that $\widetilde \phi=1$ on supp$[P,\phi]$
and supp$\,\widetilde\phi\subset\Omega_\delta$. 
On the other hand, $A$ is a finite sum of functions $a$ for which the estimate (\ref{eq:2.8}) holds. Therefore, summing up all these
estimates we arrive at 
\begin{equation}\label{eq:2.10}
\|{\rm Op}_h(A)(\phi u)\|\lesssim \lambda\|P(\phi u)\|_{-1}+h\|\phi u\|$$
$$\lesssim \lambda\|\phi Pu\|+\lambda\|[P,\phi]\widetilde \phi u\|_{-1}+h\|\phi u\|$$
$$\lesssim \lambda\|\phi Pu\|+\|\widetilde \phi u\|+h\|\phi u\|$$
$$\lesssim \lambda\|Pu\|+\|u\|_{L^2(\Omega_\delta)}+h\|\phi u\|.
\end{equation}
Clearly, the estimate (\ref{eq:2.7}) follows from (\ref{eq:2.9}) and (\ref{eq:2.10}) for $h$ small enough.
\eproof

\begin{lemma} \label{2.4} We have the estimates
\begin{equation}\label{eq:2.11}
\|\omega\|_0\lesssim \|v\|+\|f\|_{1,0}+\|u\|_1,
\end{equation}
\begin{equation}\label{eq:2.12}
h^{1/2}\|f\|_0\lesssim \varepsilon^{-1}\|u\|+\varepsilon\|u\|_1
\end{equation}
for every $0<\varepsilon\le 1$.
\end{lemma}

{\it Proof.}
Let ${\cal V}\subset\mathbb{R}^d$ be a small open domain such that ${\cal V}^0:={\cal V}\cap\Gamma\neq\emptyset$. 
Let $(x_1,x')\in {\cal V}^+:={\cal V}\cap\Omega$, $0<x_1\ll 1$, $x'=(x_2,...,x_d)\in{\cal V}^0$, be the local normal geodesic coordinates near the boundary. In these coordinates the principal symbol of the Euclidean Laplacian $-\Delta$ is equal to 
$\xi_1^2+r(x,\xi')$, where $(\xi_1,\xi')$ are the dual variables to $(x_1,x')$, and $r$ is a homogeneous polynomial of order two and satisfies $C_1|\xi'|^2\le r\le C_2|\xi'|^2$
with some constants $C_1,C_2>0$. 
Therefore, the principal symbol of the positive Laplace-Beltrami operator on $\Gamma$ is equal to
$r_0(x',\xi')=r(0,x',\xi')$. Note that $\Gamma$ can be considered as a Riemannian manifold without boundary with a 
Riemannian metric induced by the Euclidean one. 

Let ${\cal V}_1\subset{\cal V}$ be a small open domain such that ${\cal V}_1^0:={\cal V_1}\cap\Gamma\neq\emptyset$. 
Choose a function $\psi\in C_0^\infty({\cal V})$, $0\le\psi\le 1$, such that $\psi=1$ on ${\cal V}_1$. 
 We will now write the operator $-P$ in the coordinates $x=(x_1,x')$.  
 Denote ${\cal D}_{x_j}=-ih\partial_{x_j}$. We can write
$$-P=\frac{c(x)}{n(x)}\left({\cal D}_{x_1}^2+r(x,{\cal D}_{x'})\right)-1+h{\cal R}(x,{\cal D}_{x}),$$
where ${\cal R}$ is a first-order differential operator. 
 Set $u^\flat:=\psi(1-\phi)u$ and introduce the function
$$F(x_1)=\left\|{\cal D}_{x_1}u^\flat\right\|_0^2
-\left\langle r(x_1,\cdot,{\cal D}_{x'})u^\flat,u^\flat\right\rangle_0+\left\langle \widetilde n(x_1,\cdot)u^\flat,u^\flat\right\rangle_0,$$
where $\widetilde n=c^{-1}n$. Clearly, 
\begin{equation}\label{eq:2.13}
{\rm Re}\,F(0)\ge\|\psi_0\omega\|_0^2-C\|f\|_{1,0}^2, \quad C>0,
\end{equation}
where $\psi_0=\psi|_\Gamma$. 
 On the other hand,
\begin{equation}\label{eq:2.14}
F(0)=-\int_0^{\delta}F'(x_1)dx_1
\end{equation} 
where $F'$ denotes the first derivative with respect to $x_1$. We will now bound ${\rm Re}\,F(0)$ from above. 
To this end we will compute $F'(x_1)$.
We have
$$F'(x_1)=-2{\rm Re}\,\left\langle ({\cal D}_{x_1}^2+r-\widetilde n)u^\flat,
\partial_{x_1}u^\flat\right\rangle_0-\left\langle (r'-\widetilde n')u^\flat,u^\flat\right\rangle_0$$
$$=-2h^{-1}{\rm Im}\,\left\langle\widetilde n(P+h{\cal R})u^\flat,{\cal D}_{x_1}u^\flat\right\rangle_0-\left\langle (r'-\widetilde n')u^\flat,u^\flat\right\rangle_0.$$
Hence
\begin{equation}\label{eq:2.15}
|F'(x_1)|\lesssim h^{-2}\|Pu^\flat\|_0^2+\|u^\flat\|_{1,0}^2+\|{\cal D}_{x_1}u^\flat\|_0^2.
\end{equation} 
By (\ref{eq:2.14}) and (\ref{eq:2.15}) we obtain
\begin{equation}\label{eq:2.16}
{\rm Re}\,F(0)\le \int_0^{\delta}|F'(x_1)|dx_1\lesssim  h^{-2}\|Pu^\flat\|^2+\|u^\flat\|_1^2.
\end{equation}
By (\ref{eq:2.13}) and (\ref{eq:2.16}),
\begin{equation}\label{eq:2.17}
\left\|\psi_0\omega\right\|_0\lesssim h^{-1}\|Pu^\flat\|+\|u^\flat\|_1+\|f\|_{1,0}$$
$$\lesssim 
\lambda\|P((1-\phi)u)\|+\|(1-\phi)u\|_1+\|f\|_{1,0}$$
$$\lesssim 
\lambda\|Pu\|+\|u\|_1+\|f\|_{1,0}.
\end{equation}
Since $\Gamma$ is compact, there exist a finite number of smooth functions $\psi_i$, $0\le\psi_i\le 1$, $i=1,...,I,$ such that 
$1=\sum_{i=1}^I\psi_i$ and (\ref{eq:2.17}) holds with $\psi_0$ replaced by each $\psi_i$. Therefore, 
 the estimate (\ref{eq:2.11}) is obtained by summing up all such estimates (\ref{eq:2.17}).
 
 To prove (\ref{eq:2.12}) observe that 
 $$-\frac{d}{dx_1}\|u^\flat(x_1,\cdot)\|_0^2=-2{\rm Re}\langle u^\flat(x_1,\cdot),\partial_{x_1}u^\flat(x_1,\cdot)\rangle_0$$ 
 $$\le h^{-1}\varepsilon^{-2}\|u^\flat(x_1,\cdot)\|_0^2+h^{-1}\varepsilon^2\|{\cal D}_{x_1}u^\flat(x_1,\cdot)\|_0^2.$$
 Hence
 $$\|\psi_0f\|_0^2=\|u^\flat(0,\cdot)\|_0^2=-\int_0^\delta \frac{d}{dx_1}\|u^\flat(x_1,\cdot)\|_0^2dx_1$$
 $$\le h^{-1}\varepsilon^{-2}\|u^\flat\|^2+h^{-1}\varepsilon^2\|{\cal D}_{x_1}u^\flat\|^2\lesssim h^{-1}\varepsilon^{-2}\|u\|^2+h^{-1}\varepsilon^2\|u\|_1^2,$$
  which clearly implies (\ref{eq:2.12}).
\eproof

By the Green formula we have
\begin{equation}\label{eq:2.18}
\langle(\lambda^2n+i\lambda m)u-\lambda v,u\rangle=
\langle -\nabla c\nabla u,u\rangle=\int_\Omega c|\nabla u|^2+\lambda\langle f,c\omega\rangle_0.
\end{equation}
Taking the imaginary part of this identity we get
\begin{equation}\label{eq:2.19}
\langle mu,u\rangle={\rm Im}\langle v,u\rangle+{\rm Im}\langle f,c\omega\rangle_0.
\end{equation}
In view of assumption (\ref{eq:1.6}), there is a constant $m_0>0$ such that $|m|\ge m_0$ on $\Omega_{\delta}$,
provided $\delta$ is taken small enough. Thus we deduce from (\ref{eq:2.19}),
\begin{equation}\label{eq:2.20}
\|u\|_{L^2(\Omega_{\delta})}^2\lesssim \|v\|\|u\|+|{\rm Im}\,\langle f,c\omega\rangle_0|.
\end{equation}
Clearly, the estimate (\ref{eq:2.3}) follows from (\ref{eq:2.7}) and (\ref{eq:2.20}). 
On the other hand, taking the real part of (\ref{eq:2.18}) leads to the estimate
\begin{equation}\label{eq:2.21}
\|u\|_1\lesssim h\|v\|+\|u\|+h^{1/2}\|f\|_0^{1/2}\|\omega\|_0^{1/2}.
\end{equation}
By (\ref{eq:2.3}) and (\ref{eq:2.21}),
\begin{equation}\label{eq:2.22}
\|u\|_1\lesssim \|v\|+h^{1/2}\|f\|_0^{1/2}\|\omega\|_0^{1/2}+|{\rm Im}\,\langle f,c\omega\rangle_0|^{1/2}.
\end{equation}
Hence
\begin{equation}\label{eq:2.23}
\|u\|_1\lesssim \|v\|+\varepsilon^{-1}\|f\|_0+(\varepsilon+h)\|\omega\|_0
\end{equation}
for every $0<\varepsilon\le 1$. 
We now get (\ref{eq:2.2}) by combining (\ref{eq:2.11}) and (\ref{eq:2.23}), and taking $\varepsilon$ and $h$ 
small enough in order to absorb the term $(\varepsilon+h)\|\omega\|_0$ in the right-hand side. It is also easy to see that 
(\ref{eq:2.4}) follows from (\ref{eq:2.2}) and (\ref{eq:2.22}). To prove (\ref{eq:2.5}) we will use (\ref{eq:2.12}) together with 
(\ref{eq:2.3}) and (\ref{eq:2.4}).
We get
\begin{equation}\label{eq:2.24}
\|f\|_0\lesssim \varepsilon^{-1}|\lambda|^{1/2}\|v\|+\varepsilon^{-1}|\lambda|^{1/2}|{\rm Im}\,\langle f,c\omega\rangle_0|^{1/2}
+\varepsilon \|f\|_{1,0}.
\end{equation}
On the other hand, we have 
\begin{equation}\label{eq:2.25}
\|f\|_{1,0}\lesssim \|{\rm Op}_h(\chi)f\|_{1,0}+\|{\rm Op}_h(1-\chi)f\|_{1,0}
\lesssim \|f\|_0+\|{\rm Op}_h(1-\chi)f\|_{1,0}.
\end{equation}
We now combine (\ref{eq:2.24}) and (\ref{eq:2.25}) and take $\varepsilon$ small enough in order to absorb the term
$\varepsilon \|f\|_0$ 
in the right-hand side. This clearly leads to the estimate (\ref{eq:2.5}).
\eproof

Let $\chi\in C^\infty(T^*\Gamma)$ be supported in the hyperbolic region, ${\cal H}$, of the boundary value problem (\ref{eq:2.1}), that is,
$$ {\rm supp}\,\chi\subset{\cal H}:=\{(x',\xi')\in T^*\Gamma:r_0(x',\xi')<\widetilde n_0(x')\},$$
where $\widetilde n_0=\widetilde n|_\Gamma$, $\widetilde n=n/c$. 
We will now prove the following

\begin{prop} \label{2.5} Under the assumptions of Theorem 2.1, we have the estimates
\begin{equation}\label{eq:2.26}
\|{\rm Op}_h(\chi)f\|_0+\|{\rm Op}_h(\chi)\omega\|_0\lesssim
\|v\|+|{\rm Im}\,\langle f,c\omega\rangle_0|^{1/2}+h^{1/2}\|f\|_{1,0},
\end{equation}
\begin{equation}\label{eq:2.27}
\|f\|_{1,0}\lesssim
\|v\|+|{\rm Im}\,\langle f,c\omega\rangle_0|^{1/2}+\|{\rm Op}_h(1-\chi)f\|_{1,0}.
\end{equation}
\end{prop}

{\it Proof.} By (\ref{eq:2.25}) we have 
\begin{equation}\label{eq:2.28}
\|f\|_{1,0}\lesssim \|{\rm Op}_h(\chi)f\|_0+\|{\rm Op}_h(1-\chi)f\|_{1,0}.
\end{equation}
Clearly, the estimate (\ref{eq:2.27}) follows from (\ref{eq:2.26}) and (\ref{eq:2.28})
for $h$ small enough. On the other hand, 
it is easy to see that (\ref{eq:2.26}) follows from (\ref{eq:2.4}) and the following 

\begin{lemma} \label{2.6} We have the estimate
\begin{equation}\label{eq:2.29}
\|{\rm Op}_h(\chi)f\|_0+\|{\rm Op}_h(\chi)\omega\|_0\lesssim\|v\|+\|u\|_1+h\|f\|_0.
\end{equation}
\end{lemma}

{\it Proof.} Let $u^\flat$ be as in the proof of Lemma 2.4 and set 
$$u^\sharp={\rm Op}_h(\chi)u^\flat={\rm Op}_h(\chi)\psi(1-\phi)u.$$ 
Clearly,
\begin{equation}\label{eq:2.30}
\|Pu^\sharp\|\lesssim \|P((1-\phi)u)\|+\|[P,{\rm Op}_h(\chi)\psi](1-\phi)u\|$$ 
$$\lesssim \|P((1-\phi)u)\|+h\|(1-\phi)u\|_1$$
$$\lesssim \|Pu\|+h\|u\|_1\lesssim h\|v\|+h\|u\|_1.
\end{equation}
We define the function $F(x_1)$ as in the proof of Lemma 2.4 replacing $u^\flat$ by $u^\sharp$. Observe now that the choice of
$\chi$ guarantees that
$$\widetilde n(0,x')-r_0(x',\xi')\ge C>0$$
on supp$\,\chi$. Therefore, by G{\aa}rding's inequality we have 
\begin{equation}\label{eq:2.31}
{\rm Re}\,\left\langle (\widetilde n(0,\cdot)- r_0(\cdot,{\cal D}_{x'})){\rm Op}_h(\chi)\psi_0f,{\rm Op}_h(\chi)\psi_0f\right\rangle_0
\ge C_1\|{\rm Op}_h(\chi)\psi_0f\|_0^2,\quad C_1>0.
\end{equation}
Since
$${\cal D}_{x_1}u^\sharp|_{x_1=0}=-i{\rm Op}_h(\chi)\psi_0\omega-ih{\rm Op}_h(\chi)\psi_1f,$$
where $\psi_0=\psi|_{x_1=0}$, $\psi_1=\partial_{x_1}\psi|_{x_1=0}$, we deduce from (\ref{eq:2.31}),
\begin{equation}\label{eq:2.32}
{\rm Re}\,F(0)\ge C_1\|{\rm Op}_h(\chi)\psi_0f\|_0^2+\|{\rm Op}_h(\chi)\psi_0\omega\|_0^2-{\cal O}(h^2)\|f\|_0^2$$ 
$$\ge C_1\|\psi_0{\rm Op}_h(\chi)f\|_0^2+\|\psi_0{\rm Op}_h(\chi)\omega\|_0^2-{\cal O}(h^2)\|f\|_0^2.
\end{equation}
On the other hand, the upper bound (\ref{eq:2.16}) still holds with $u^\flat$ replaced by $u^\sharp$.
This fact together with (\ref{eq:2.30}) and (\ref{eq:2.32}) imply (\ref{eq:2.29}).
\eproof

Given a parameter $0<\varepsilon\ll 1$, independent of $\lambda$, choose a function  $\chi_\varepsilon\in C_0^\infty(T^*\Gamma)$
such that $\chi_\varepsilon=1$ in the region $\{|cr_0/n-1|\le\varepsilon\}$ and 
$\chi_\varepsilon=0$ in $T^*\Gamma\setminus\{|cr_0/n-1|\le 2\varepsilon\}$.
In the same way as above we get the following

\begin{prop} \label{2.7} Under the assumptions of Theorem 2.1, we have the estimate
\begin{equation}\label{eq:2.33}
\|{\rm Op}_h(\chi_\varepsilon)\omega\|_0\lesssim 
\|v\|+|{\rm Im}\,\langle f,c\omega\rangle_0|^{1/2}+(\varepsilon+h^{1/2})\|f\|_{1,0}.
\end{equation}
\end{prop}

{\it Proof.} It is easy to see that the proposition follows from (\ref{eq:2.4}) and the following 

\begin{lemma} \label{2.8} We have the estimate
\begin{equation}\label{eq:2.34}
\|{\rm Op}_h(\chi_\varepsilon)\omega\|_0\lesssim\|v\|+\|u\|_1+(\varepsilon+h)\|f\|_0.
\end{equation}
\end{lemma}

{\it Proof.} We will proceed in the same way as in the proof of Lemma 2.6 with $\chi$ replaced by $\chi_\varepsilon$ making the following modification. Since in this case the function $\chi_\varepsilon$ is no longer supported in the hyperbolic region, we do not have
the G{\aa}rding inequality (\ref{eq:2.31}) fulfilled anymore. Instead, since $\widetilde n(0,x')-r_0(x',\xi')={\cal O}(\varepsilon)$ on 
supp$\,\chi_\varepsilon$, we have the bound
$$\left\|(\widetilde n(0,\cdot)- r_0(\cdot,{\cal D}_{x'})){\rm Op}_h(\chi_\varepsilon)\psi_0f\right\|_0
\lesssim (\varepsilon+h)\|\psi_0f\|_0.$$
Therefore in this case the function ${\rm Re}\,F(0)$ is lower bounded as follows
\begin{equation}\label{eq:2.35}
{\rm Re}\,F(0)\ge \|{\rm Op}_h(\chi_\varepsilon)\psi_0\omega\|_0^2-{\cal O}((\varepsilon+h)^2)\|f\|_0^2$$ 
$$\ge \|\psi_0{\rm Op}_h(\chi_\varepsilon)\omega\|_0^2-{\cal O}((\varepsilon+h)^2)\|f\|_0^2.
\end{equation}
The estimate (\ref{eq:2.34}) follows from (\ref{eq:2.35}) and the upper bound (\ref{eq:2.16}) adapted to our case.
\eproof

\section{A priori estimates when $m\equiv 0$}

In this section we consider the equation
\begin{equation}\label{eq:3.1}
\left\{
\begin{array}{l}
(\nabla c(x)\nabla+\lambda^2n(x))u=\lambda v\quad \mbox{in}\quad\Omega,\\
u=f\quad\mbox{on}\quad\Gamma,\\
\end{array}
\right.
\end{equation}
where $\lambda\in\Lambda_k$ and $c,n\in C^\infty(\overline\Omega)$ are real-valued functions satisfying $c(x)>0$, $n(x)>0$ for all $x\in \overline\Omega$.  We keep the same notations as in the previous section. We have the following

\begin{Theorem} \label{3.1} Let $\Gamma$ be $g-$strictly concave and suppose that the condition (\ref{eq:1.7})
is fulfilled with $g_j$ replaced by $g=\frac{c(x)}{n(x)}|\xi|^2$. Let $u\in H^2(\Omega)$ satisfy equation (\ref{eq:3.1}) and set 
 $\omega=h\partial_\nu u|_\Gamma$. Then there are constants $k,\lambda_0>0$ such that for all $\lambda\in\Lambda_k$, 
 ${\rm Re}\,\lambda\ge\lambda_0$, we have the estimate
\begin{equation}\label{eq:3.2}
\|u\|_1\lesssim \|v\|+\|f\|_0+\|\omega\|_0.
\end{equation}
\end{Theorem}

In the same way as at the begining of the proof of Theorem 2.1 one can see that the estimate (\ref{eq:3.2}) for complex
$\lambda$ follows from (\ref{eq:3.2}) for real $\lambda\ge\lambda_0$. In other words, it suffices to prove 
Theorem 3.1 for real $\lambda\gg 1$, only. Let $\phi_2\in C_0^\infty(\mathbb{R}^d)$ be independent of $\lambda$ and such that $\phi_2=1$ in 
$\Omega\setminus\Omega_{3\delta}$ and supp$\,\phi_2\subset\Omega\setminus\Omega_{2\delta}$. 
The estimate (\ref{eq:3.2}) for real $\lambda\gg 1$ follows from Proposition 2.2 (which clearly holds when $m\equiv 0$) and the following

\begin{prop} \label{3.2} Let $\Gamma$ be $g-$strictly concave. Then we have the estimate
\begin{equation}\label{eq:3.3}
\|(1-\phi_2)u\|\lesssim \|v\|+\|f\|_0+\|\omega\|_0+h^{1/2}\|u\|_1,
\end{equation}
provided the parameter $\delta$ is taken small enough.
\end{prop}

Indeed, combining the estimates (\ref{eq:2.7}) and (\ref{eq:3.3}) leads to
\begin{equation}\label{eq:3.4}
\|u\|\lesssim \|v\|+\|f\|_0+\|\omega\|_0+h^{1/2}\|u\|_1. 
\end{equation}
On the other hand, clearly the estimate (\ref{eq:2.21}) still holds when $m\equiv 0$, so we have
\begin{equation}\label{eq:3.5}
\|u\|_1\lesssim \|u\|+\|v\|+\|f\|_0+\|\omega\|_0.
\end{equation}
Therefore, 
the estimate (\ref{eq:3.2}) follows from (\ref{eq:3.4}) and (\ref{eq:3.5}) after absorbing the term $h^{1/2}\|u\|_1$
by taking $h$ small enough. 

Note that the above proposition is in fact Proposition 2.2 of \cite{kn:CPV} and we refer the reader to Section 2 of \cite{kn:CPV} for the proof.
In the next sections we will also need the following

\begin{Theorem} \label{3.3} Let $u\in H^2(\Omega)$ satisfy equation (\ref{eq:3.1}). 
Then the function $\omega=h\partial_\nu u|_\Gamma$ satisfies the estimate
\begin{equation}\label{eq:3.6}
\|\omega\|_0\lesssim |{\rm Im}\,\lambda|^{-1}\left(\|v\|+\|f\|_{1,0}\right),\quad {\rm Im}\,\lambda\neq 0.
\end{equation}
\end{Theorem}

 The theorem follows from Lemma 2.4 (which clearly holds when $m\equiv 0$) and the following  

\begin{lemma} \label{3.4} We have the estimate
\begin{equation}\label{eq:3.7}
\|u\|_1\lesssim |{\rm Im}\,\lambda|^{-1}\|v\|+|{\rm Im}\,\lambda|^{-1/2}\|f\|_{0}^{1/2}\|\omega\|_0^{1/2}.
\end{equation}
\end{lemma}

This lemma follows easily from the Green formula (see Lemma 3.1 of \cite{kn:V2}). 

\section{The Dirichlet-to-Neumann map}

Let $u$ solve equation (\ref{eq:2.1}) with $v\equiv 0$ and define the Dirichlet-to-Neumann map
$${\cal N}(\lambda):H^1(\Gamma)\to L^2(\Gamma)$$
by
$${\cal N}(\lambda)f:=h\partial_\nu u|_\Gamma.$$
Clearly, under the conditions of Theorem 2.1, by (\ref{eq:2.2}) we have that the Dirichlet-to-Neumann map in this case satisfies the estimate 
\begin{equation}\label{eq:4.1}
\|{\cal N}(\lambda)f\|_0\lesssim \|f\|_{1,0}.
\end{equation}
On the other hand, when $m\equiv 0$, by (\ref{eq:3.6}) we have
the estimate
\begin{equation}\label{eq:4.2}
\|{\cal N}(\lambda)f\|_0\lesssim |{\rm Im}\,\lambda|^{-1}\|f\|_{1,0},\quad {\rm Im}\,\lambda\neq 0.
\end{equation}
Let $\chi,\eta\in C^\infty(T^*\Gamma)$ be compactly supported functions such that $\eta=1$ on supp$\,\chi$. 
 In Section 6 we will need the following

\begin{lemma} \label{4.1} Under the assumptions of Theorem 2.1, we have the estimates
\begin{equation}\label{eq:4.3}
\|[{\cal N}(\lambda),{\rm Op}_h(\chi)]f\|_0\lesssim h^{1/2}\|f\|_{1,0}
+|{\rm Im}\,\langle c{\cal N}(\lambda)f,f\rangle_0|^{1/2},
\end{equation}
\begin{equation}\label{eq:4.4}
\|{\rm Op}_h(1-\eta){\cal N}(\lambda){\rm Op}_h(\chi)f\|_0\lesssim h^{1/2}\|f\|_{1,0}
+|{\rm Im}\,\langle c{\cal N}(\lambda)f,f\rangle_0|^{1/2}.
\end{equation}
\end{lemma}

{\it Proof.} Let $u$ solve equation (\ref{eq:2.1}) with $v\equiv 0$. Then the function
$\widetilde u={\rm Op}_h(\chi)(1-\phi)u$ solves equation (\ref{eq:2.1}) with $v$ and $f$ replaced by $\widetilde v$
and $\widetilde f$, respectively, where
$$\widetilde v=\lambda^{-1}[\nabla c(x)\nabla+\lambda^2n(x)+i\lambda m(x),{\rm Op}_h(\chi)(1-\phi)]u,$$
$$\widetilde f={\rm Op}_h(\chi)f.$$
Let $w$ solve equation (\ref{eq:2.1}) with $v\equiv 0$ and $f$ replaced by $\widetilde f$.
Then the function $\widetilde u-w$ solves equation (\ref{eq:2.1}) with $v$ replaced by $\widetilde v$ and $f=0$.
Moreover, we have
$$h\partial_\nu(\widetilde u-w)|_\Gamma=-[{\cal N}(\lambda),{\rm Op}_h(\chi)]f.$$
Therefore, the estimate (\ref{eq:2.2}) in this case leads to 
$$\|[{\cal N}(\lambda),{\rm Op}_h(\chi)]f\|_0\lesssim \|\widetilde v\|\lesssim \|u\|_1.$$
Hence the estimate (\ref{eq:4.3}) follows from (\ref{eq:2.4}). To prove (\ref{eq:4.4}) we will use that
$${\rm Op}_h(1-\eta){\rm Op}_h(\chi)={\cal O}(h^\infty):L^2(\Gamma)\to L^2(\Gamma).$$
Thus, in view of (\ref{eq:4.1}), we get
$$\|{\rm Op}_h(1-\eta){\cal N}(\lambda){\rm Op}_h(\chi)f\|_0\lesssim \|{\rm Op}_h(1-\eta)[{\cal N}(\lambda),{\rm Op}_h(\chi)]f\|_0$$
$$+\|{\rm Op}_h(1-\eta){\rm Op}_h(\chi){\cal N}(\lambda)f\|_0$$
$$\lesssim \|[{\cal N}(\lambda),{\rm Op}_h(\chi)]f\|_0+h^\infty\|{\cal N}(\lambda)f\|_0$$
$$\lesssim \|[{\cal N}(\lambda),{\rm Op}_h(\chi)]f\|_0+h^\infty\|f\|_{1,0}.$$
We now obtain (\ref{eq:4.4}) from (\ref{eq:4.3}).
\eproof

Denote by ${\cal N}(\lambda)^*$ the adjoint of ${\cal N}(\lambda)$ with respect to the scalar product
$\langle\cdot,\cdot\rangle_0$ in $L^2(\Gamma)$. 
In Section 5 we will need the following
\begin{lemma} \label{4.2} When $m\equiv 0$ the Dirichlet-to-Neumann map satisfies the identity
\begin{equation}\label{eq:4.5}
{\cal N}(\lambda)^*c_0=c_0{\cal N}(\overline\lambda),
\end{equation}
where $c_0=c|_\Gamma$. 
\end{lemma}

{\it Proof.} Given any $f_1,f_2\in L^2(\Gamma)$, let $u_1$ be the solution of equation (\ref{eq:3.1})
with $v\equiv 0$ and $f$ replaced by $f_1$, and let $u_2$ be the solution of equation (\ref{eq:3.1})
with $v\equiv 0$, $f$ replaced by $f_2$ and $\lambda$ replaced by $\overline\lambda$. By the Green formula we have 
$$0=-\langle\nabla c\nabla u_1,u_2\rangle+\langle u_1,\nabla c\nabla u_2\rangle
=\langle c_0\partial_\nu u_1|_\Gamma,f_2\rangle_0-\langle f_1,c_0\partial_\nu u_2|_\Gamma\rangle_0$$
 $$=h^{-1}\langle c_0{\cal N}(\lambda)f_1,f_2\rangle_0-h^{-1}\langle f_1,c_0{\cal N}(\overline\lambda)f_2\rangle_0,$$
which clearly implies (\ref{eq:4.5}).
\eproof

\section{Parametrix of the Dirichlet-to-Neumann map in the elliptic region}

Let $\eta\in C^\infty(T^*\Gamma)$ be such that $1-\eta$ is supported in the elliptic region, ${\cal E}$, of the boundary value problem (\ref{eq:2.1}), that is,
$${\rm supp}(1-\eta)\subset{\cal E}:=\{(x',\xi')\in T^*\Gamma:r_0(x',\xi')>\widetilde n_0(x')\}.$$
For $(x',\xi')\in {\cal E}$ set 
$$\rho(x',\xi',z)=\sqrt{r_0(x',\xi')-z\widetilde n_0(x')},\quad {\rm Re}\,\rho>0,$$
where $z=(h\lambda)^2=1-\theta^2+2i\theta$, $\theta=h{\rm Im}\,\lambda$. In fact, on ${\rm supp}(1-\eta)$ we have the lower bound
\begin{equation}\label{eq:5.1}
{\rm Re}\,\rho\ge C\langle\xi'\rangle, \quad C>0.
\end{equation} 
Given any $N>1$, set ${\cal L}_N:=\{\lambda\in\mathbb{C}, |{\rm Im}\,\lambda|\le |\lambda|^{-N},\,{\rm Re}\,\lambda\ge 1\}$.
Our goal in this section is to prove the following

\begin{Theorem} Under the conditions of Theorem 2.1 
we have the estimate
\begin{equation}\label{eq:5.2}
\left\|{\cal N}(\lambda){\rm Op}_h(1-\eta)f+{\rm Op}_h(\rho(1-\eta))f\right\|_0\lesssim h\|f\|_0 
\end{equation}
for $\lambda\in\Lambda_k$. When $m\equiv 0$ the estimate (\ref{eq:5.2}) still holds for $\lambda\in\Lambda_k\setminus{\cal L}_N$
without any conditions.
 \end{Theorem}

{\it Proof.} The theorem follows from the parametrix construction carried out in \cite{kn:V1}. In what follows we will recall it.
In fact, in \cite{kn:V1} the case $m\equiv 0$ is considered, but it is easy to see that the presence of the function
$m$ does not change anything. Indeed, the eikonal equation does not depend on $m$ and only the transport equations do.
Note also that it suffices to build the parametrix locally and then sum up all pieces.

  Let $(x_1,x')\in {\cal V}^+$ be the local normal geodesic coordinates near the boundary. 
Take a function $\chi\in C^\infty(T^*\Gamma)$, $0\le\chi\le 1$, such that $\pi_{x'}({\rm supp}\,\chi)\subset {\cal V}^0$, where 
$\pi_{x'}:T^*\Gamma\to\Gamma$ denotes the projection $(x',\xi')\to x'$. Moreover, we require that 
 $\chi\in S^0(\Gamma)$ with ${\rm supp}\,\chi\subset{\rm supp}(1-\eta)$. 
We will be looking for a parametrix of the solution to 
equation (\ref{eq:2.1}) (with $v\equiv 0$) in the form
$$\widetilde u=\phi_0(x_1)(2\pi h)^{-d+1}\int\int e^{\frac{i}{h}(\langle y',\xi'\rangle+\varphi(x,\xi',z))}a(x,\xi',z,h)f(y')d\xi'dy',$$
where $\phi_0\in C_0^\infty(\mathbb{R})$,
$\phi_0(t)=1$ for $|t|\le \delta/2$, $\phi_0(t)=0$ for $|t|\ge\delta$. Here $0<\delta\ll 1$ is a small parameter independent of $\lambda$.   
 We require that $\widetilde u$ satisfies the boundary condition
$\widetilde u={\rm Op}_h(\chi)f$ on $x_1=0$. The phase $\varphi$ and the amplitude $a$ are choosen in such a way that the function
$\widetilde u$ satisfies equation (\ref{eq:2.1}) mod ${\cal O}(h^M)$, where $M\ge 1$ is an arbitrary integer. 
The phase function satisfies 
$$\varphi|_{x_1=0}=-\langle x',\xi'\rangle$$
as well as the eikonal equation
\begin{equation}\label{eq:5.3}
(\partial_{x_1}\varphi)^2+r(x,\nabla_{x'}\varphi)-z\widetilde n(x)=x_1^M\Psi_M,
\end{equation}
where the function $\left|\Psi_M\right|$ is bounded as $x_1\to 0$. It is shown in Section 4 of \cite{kn:V1} that 
(\ref{eq:5.3}) has a solution of the form
$$\varphi=\sum_{j=0}^Mx_1^j\varphi_j,$$
where the functions $\varphi_j$ do not depend on $x_1$, $\varphi_0=-\langle x',\xi'\rangle$,
$\varphi_1=i\rho$. It follows from (\ref{eq:5.1}) that
\begin{equation}\label{eq:5.4}
{\rm Im}\,\varphi\ge Cx_1\langle\xi'\rangle/2,
\end{equation} 
for $0\le x_1\le\delta$, provided $\delta$ is taken small enough. 
The amplitude is of the form
$$a=\sum_{j=0}^Mh^ja_j,$$
where the functions $a_j$ do not depend on $h$, $a_0|_{x_1=0}=\chi$. Then all functions $a_j$ can be determined from the transport equations
and we have $a_j\in S^{-j}(\Gamma)$ uniformly in $x_1$ and $z$ (see Section 4 of \cite{kn:V1}). Clearly, we have
\begin{equation}\label{eq:5.5}
h\partial_\nu\widetilde u|_{x_1=0}={\rm Op}_h(b_M)f,
\end{equation}
where
$$b_M=ia\frac{\partial\varphi}{\partial x_1}|_{x_1=0}+h\frac{\partial a}{\partial x_1}|_{x_1=0}=-\chi\rho+h\sum_{j=0}^Mh^j\frac{\partial a_j}{\partial x_1}|_{x_1=0}.$$
Hence $h^{-1}(b_M+\chi\rho)\in S^0(\Gamma)$ uniformly in $h$. This implies
\begin{equation}\label{eq:5.6}
{\rm Op}_h(b_M+\chi\rho)={\cal O}(h):L^2(\Gamma)\to L^2(\Gamma).
\end{equation}
On the other hand, the function
$$\widetilde v=(\nabla c(x)\nabla+\lambda^2n(x)+i\lambda m(x))\widetilde u$$
is of the form
$$\widetilde v=(2\pi h)^{-d+1}\int\int e^{\frac{i}{h}(\langle y',\xi'\rangle+\varphi(x,\xi',z))}V_M(x,\xi',z,h)f(y')d\xi'dy',$$
where $V_M=V_M^{(1)}+\phi_0(x_1)V_M^{(2)}$, 
$$V_M^{(1)}=[\nabla c(x)\nabla,\phi_0(x_1)]a,$$
$$V_M^{(2)}=e^{-i\varphi/h}(\nabla c(x)\nabla+\lambda^2n(x)+i\lambda m(x))e^{i\varphi/h}a.$$
As shown in Section 4 of \cite{kn:V1}, the functions $a_j$ can be choosen in such a way that the function $V_M^{(2)}$ is of the form
\begin{equation}\label{eq:5.7}
V_M^{(2)}=x_1^MA_M+h^MB_M,
\end{equation}
where $A_M$ and $B_M$ are smooth functions. More precisely, since $\chi$ is supported in the elliptic region, we have 
$A_M\in S^2(\Gamma)$, $B_M\in S^{1-M}(\Gamma)$ uniformly in $h$, $z$ and $0<x_1\le\delta$
(see Proposition 3.4 of \cite{kn:V1}). Note that in view of (\ref{eq:5.4}) we have the bound
$$\left|x_1^Me^{i\varphi/h}\right|\lesssim h^M\langle\xi'\rangle^{-M}.$$
Thus we get that the function $V_M$ satisfies the bound
\begin{equation}\label{eq:5.8}
\left|V_Me^{i\varphi/h}\right|\lesssim h^M\langle\xi'\rangle^{-M+1}.
\end{equation}
By (\ref{eq:5.8}) we obtain the estimate
\begin{equation}\label{eq:5.9}
\|\widetilde v\|\le C_Mh^{M/2}\|f\|_0,
\end{equation}
provided $M$ is taken big enough. Let $u$ solve equation (\ref{eq:2.1}) with $v\equiv 0$ and $f$ replaced by ${\rm Op}_h(\chi)f$. Then the function
$\widetilde u-u$ solves equation (\ref{eq:2.1}) with $v$ replaced by $\widetilde v$ and $f=0$. Therefore, under the conditions of
Theorem 2.1, by (\ref{eq:2.2}) together with (\ref{eq:5.5}) and (\ref{eq:5.9}) we get
\begin{equation}\label{eq:5.10}
\left\|{\cal N}(\lambda){\rm Op}_h(\chi)f-{\rm Op}_h(b_M)f\right\|_0=
\|h\partial_{x_1}(\widetilde u-u)|_{x_1=0}\|\lesssim\|\widetilde v\|\lesssim h^{M/2}\|f\|_0.
\end{equation}
By (\ref{eq:5.6}) and (\ref{eq:5.10}),
\begin{equation}\label{eq:5.11}
\left\|{\cal N}(\lambda){\rm Op}_h(\chi)f+{\rm Op}_h(\chi\rho)f\right\|_0\lesssim h\|f\|_0,
\end{equation}
which implies (\ref{eq:5.2}) in this case since $1-\eta$ can be written as a finite sum of functions $\chi$ for which
(\ref{eq:5.11}) holds. Consider now the case when $m\equiv 0$. We proceed similarly with the difference that 
we use the estimate (\ref{eq:3.6}) instead of (\ref{eq:2.2}). For $\lambda\in\Lambda_k\setminus{\cal L}_N$, we obtain
\begin{equation}\label{eq:5.12}
\left\|{\cal N}(\lambda){\rm Op}_h(\chi)f-{\rm Op}_h(b_M)f\right\|_0\lesssim |{\rm Im}\,\lambda|^{-1}\|\widetilde v\|
\lesssim h^{M/2-N}\|f\|_0\lesssim h\|f\|_0,
\end{equation}
provided we take $M\ge 2N+1$. Thus we conclude that the estimate (\ref{eq:5.11}) (and hence (\ref{eq:5.2}))
 still holds in this case as long as
$\lambda\in\Lambda_k\setminus{\cal L}_N$ with a constant in front of the term in the right-hand side depending on $N$.
\eproof

Let $\chi,\eta\in C^\infty(T^*\Gamma)$ be compactly supported functions such that $\eta=1$ on supp$\,\chi$
and supp$(1-\eta)\subset{\cal E}$. We will use Theorem 5.1 to prove the following

\begin{lemma} \label{5.2} When $m\equiv 0$ we have the estimates
\begin{equation}\label{eq:5.13}
\|{\rm Op}_h(\chi){\cal N}(\lambda){\rm Op}_h(1-\eta)f\|_0\lesssim h\|f\|_0,
\end{equation}
\begin{equation}\label{eq:5.14}
\|{\rm Op}_h(1-\eta){\cal N}(\lambda){\rm Op}_h(\chi)f\|_0\lesssim h\|f\|_0,
\end{equation}
for $\lambda\in\Lambda_k\setminus{\cal L}_N$.
\end{lemma}

{\it Proof.} Since
$${\rm Op}_h(\chi){\rm Op}_h(\rho(1-\eta))={\cal O}(h^\infty):L^2(\Gamma)\to L^2(\Gamma),$$
the estimate (\ref{eq:5.13}) follows from (\ref{eq:5.2}). In view of Lemma 4.2 the adjoint of the operator
$${\cal A}:={\rm Op}_h(1-\eta){\cal N}(\lambda){\rm Op}_h(\chi)$$
is
$${\cal A}^*={\rm Op}_h(\chi)^*c_0{\cal N}(\overline\lambda)c_0^{-1}{\rm Op}_h(1-\eta)^*.$$
Choose compactly supported functions $\chi_1,\eta_1\in C^\infty(T^*\Gamma)$ such that $\eta_1=1$ on supp$\,\chi_1$, 
 supp$(1-\eta_1)\subset{\cal E}$, $\chi_1=1$ on supp$\,\chi$ and $\eta=1$ on supp$\,\eta_1$.
 The standard $h-\Psi$DO calculus give
 $${\rm Op}_h(\chi)^*c_0{\rm Op}_h(1-\chi_1)={\cal O}(h^\infty):H^{-1}(\Gamma)\to L^2(\Gamma),$$
 $${\rm Op}_h(\eta_1)c_0^{-1}{\rm Op}_h(1-\eta)^*={\cal O}(h^\infty):L^2(\Gamma)\to H^1(\Gamma).$$
 We now apply the estimate (\ref{eq:5.13}) with $\chi$, $\eta$, $\lambda$ replaced by 
 $\chi_1$, $\eta_1$, $\overline\lambda$, respectively. We will also use (\ref{eq:5.2}) with $\eta$, $\lambda$ and $\rho$ replaced by
 $\eta_1$, $\overline\lambda$ and $\overline\rho$, respectively. Note that $\overline\rho(1-\eta_1)\in S^1(\Gamma)$. 
 Thus, in view of (\ref{eq:4.2}), we get
 $$\|{\cal A}^*f\|_0\lesssim \|{\rm Op}_h(\chi_1){\cal N}(\overline\lambda){\rm Op}_h(1-\eta_1)c_0^{-1}{\rm Op}_h(1-\eta)^*f\|_0$$
 $$+\|{\cal N}(\overline\lambda){\rm Op}_h(\eta_1)c_0^{-1}{\rm Op}_h(1-\eta)^*f\|_0$$
 $$+\|{\rm Op}_h(\chi)^*c_0{\rm Op}_h(1-\chi_1){\cal N}(\overline\lambda)c_0^{-1}{\rm Op}_h(1-\eta)^*f\|_0$$
 $$\lesssim h\|c_0^{-1}{\rm Op}_h(1-\eta)^*f\|_0$$
 $$+|{\rm Im}\,\lambda|^{-1}\|{\rm Op}_h(\eta_1)c_0^{-1}{\rm Op}_h(1-\eta)^*f\|_{1,0}$$
 $$+h^\infty\|{\cal N}(\overline\lambda){\rm Op}_h(\eta_1)c_0^{-1}{\rm Op}_h(1-\eta)^*f\|_0$$
 $$+h^\infty\|{\cal N}(\overline\lambda){\rm Op}_h(1-\eta_1)c_0^{-1}{\rm Op}_h(1-\eta)^*f\|_{-1,0}$$
 $$\lesssim h\|f\|_0+h^\infty|{\rm Im}\,\lambda|^{-1}\|f\|_0$$
 $$+h^\infty\|{\rm Op}_h(\overline\rho(1-\eta_1))c_0^{-1}{\rm Op}_h(1-\eta)^*f\|_{-1,0}$$
 $$+h^\infty\|({\cal N}(\overline\lambda){\rm Op}_h(1-\eta_1)+{\rm Op}_h(\overline\rho(1-\eta_1)))c_0^{-1}{\rm Op}_h(1-\eta)^*f\|_0$$
 $$\lesssim h\|f\|_0+h^\infty|{\rm Im}\,\lambda|^{-1}\|f\|_0\lesssim h\|f\|_0,$$
 where $\|\cdot\|_{-1,0}$ denotes the semiclassical norm in $H^{-1}(\Gamma)$ and we have used that $|{\rm Im}\,\lambda|\ge h^N$.
 In other words,
 $${\cal A}^*={\cal O}(h):L^2(\Gamma)\to L^2(\Gamma),$$
 and hence so is the operator ${\cal A}$.
 \eproof
  
\section{Proof of Theorem 1.1}

Let $(u_1,u_2)$ be the solution to equation (\ref{eq:1.2}) and set $f=u_1|_\Gamma=u_2|_\Gamma$. Then we can express the restrictions of the normal derivative of $u_1$ and $u_2$ in terms of the corresponding Dirichlet-to-Neumann maps, that is,
$$h\partial_\nu u_j|_\Gamma={\cal N}_j(\lambda)f,\quad j=1,2.$$
Therefore, $\lambda$ is a transmission eigenvalue if $T(\lambda)f\equiv 0$, where
$$T(\lambda)=c_1{\cal N}_1(\lambda)-c_2{\cal N}_2(\lambda).$$
We have to show that if $\lambda$ belongs to the eigenvalue-free regions of Theorem 1.1 and $T(\lambda)f\equiv 0$,
then $f\equiv 0$. We will first prove the following

\begin{lemma} There exists a constant $k>0$ such that for $\lambda\in\Lambda_k$ we have the estimate
\begin{equation}\label{eq:6.1}
|{\rm Im}\,\langle c_1{\cal N}_1(\lambda)f,f\rangle_0|\lesssim |{\rm Im}\,\lambda|\|f\|_{1,0}^2.
\end{equation}
\end{lemma}

{\it Proof.} The Green formula applied to the second equation in (\ref{eq:1.2}) gives the identity
$$2{\rm Im}\,\lambda\langle n_2u_2,u_2\rangle={\rm Im}\,\langle c_2{\cal N}_2(\lambda)f,f\rangle_0={\rm Im}\,\langle c_1{\cal N}_1(\lambda)f,f\rangle_0.$$
On the other hand, by (\ref{eq:3.2}) and (\ref{eq:4.1}) we have
$$\|u_2\|\lesssim \|f\|_0+\|{\cal N}_2(\lambda)f\|_0\lesssim \|f\|_0+\|{\cal N}_1(\lambda)f\|_0\lesssim\|f\|_{1,0},$$
which clearly implies (\ref{eq:6.1}).
\eproof

Let $\chi\in C^\infty(T^*\Gamma)$ be of compact support such that $1-\chi$ is supported in the region $\{r_0\ge\sigma\}$, where $\sigma\gg 1$ is a constant to be fixed in the next lemma. 

\begin{lemma} For a suitable choice of $\sigma$ we have the estimate
\begin{equation}\label{eq:6.2}
\|{\rm Op}_h(1-\chi)f\|_{1,0}\lesssim (|{\rm Im}\,\lambda|^{1/2}+h^{1/2})\|f\|_{1,0}
\end{equation}
for $\lambda\in\Lambda_k\setminus{\cal L}_N$.
\end{lemma}

{\it Proof.} Choose a compactly supported function $\eta\in C^\infty(T^*\Gamma)$ such that $1-\eta$ is supported in the region $\{r_0\ge\sigma\}$
and $\chi=1$ on supp$\,\eta$. Define $\rho_j$, ${\cal H}_j$, ${\cal E}_j$, $j=1,2$, by replacing in the definition of $\rho$, 
${\cal H}$, ${\cal E}$
in Sections 4 and 5 the functions $c,n$ by $c_j,n_j$. Clearly, taking $\sigma$ big enough we can arrange that the functions 
$1-\chi$ and $1-\eta$ are supported in both elliptic regions ${\cal E}_1$ and ${\cal E}_2$. 
Then the fact that $Tf=0$ implies the identity
$${\rm Op}_h((c_1\rho_1-c_2\rho_2)(1-\chi))f$$
$$+{\rm Op}_h(1-\chi){\rm Op}_h((c_1\rho_1-c_2\rho_2)(1-\eta))f-{\rm Op}_h((c_1\rho_1-c_2\rho_2)(1-\chi))f$$
$$={\rm Op}_h(1-\chi){\rm Op}_h((c_1\rho_1-c_2\rho_2)(1-\eta))f$$
$$={\rm Op}_h(1-\chi)\left(T{\rm Op}_h(1-\eta)+{\rm Op}_h((c_1\rho_1-c_2\rho_2)(1-\eta))\right)f$$
$$+{\rm Op}_h(1-\chi)T{\rm Op}_h(\eta)f.$$
Since $(c_1\rho_1-c_2\rho_2)(1-\chi)$ and $(c_1\rho_1-c_2\rho_2)(1-\eta)$ belong to $S^1(\Gamma)$, the $h-\Psi$DO calculus give
$${\rm Op}_h(1-\chi){\rm Op}_h((c_1\rho_1-c_2\rho_2)(1-\eta))f-{\rm Op}_h((c_1\rho_1-c_2\rho_2)(1-\chi))
={\cal O}(h):H^1(\Gamma)\to L^2(\Gamma).$$
Therefore, using Theorem 5.1 together with the estimates (\ref{eq:4.4}), (\ref{eq:6.1}) and (\ref{eq:5.14}) we obtain
\begin{equation}\label{eq:6.3}
\|{\rm Op}_h((c_1\rho_1-c_2\rho_2)(1-\chi))f\|_0\lesssim (|{\rm Im}\,\lambda|^{1/2}+h^{1/2})\|f\|_{1,0}.
\end{equation}
On the other hand, we have
$$(c_1\rho_1+c_2\rho_2)(c_1\rho_1-c_2\rho_2)=c_1^2\rho_1^2-c_2^2\rho_2^2=(c_1^2-c_2^2)r_0-z(c_1n_1-c_2n_2).$$
Hence
$$\left|(c_1^2-c_2^2)r_0-z(c_1n_1-c_2n_2)\right|\lesssim\langle\xi'\rangle\left|c_1\rho_1-c_2\rho_2\right|.$$
On the other hand, in view of assumption (\ref{eq:1.3}) we have $|c_1^2-c_2^2|\ge C_0$ with some constant
$C_0>0$. Therefore, taking $\sigma$ big enough we can arrange that
$$\left|(c_1^2-c_2^2)r_0-z(c_1n_1-c_2n_2)\right|\ge C_1\langle\xi'\rangle^2,\quad C_1>0,$$
on supp$(1-\eta)$. Hence
\begin{equation}\label{eq:6.4}
\left|c_1\rho_1-c_2\rho_2\right|\ge C_2\langle\xi'\rangle,\quad C_2>0,
\end{equation}
on supp$(1-\eta)$, which implies
\begin{equation}\label{eq:6.5}
\|{\rm Op}_h(1-\chi)f\|_{1,0}\lesssim 
\|{\rm Op}_h((c_1\rho_1-c_2\rho_2)(1-\chi))f\|_0+h\|f\|_0.
\end{equation}
Clearly, the estimate (\ref{eq:6.2}) follows from (\ref{eq:6.3}) and (\ref{eq:6.5}).
\eproof

By (\ref{eq:2.5}), (\ref{eq:6.1}) and (\ref{eq:6.2}) we obtain 
$$\|f\|_{1,0}\lesssim |\lambda|^{1/2}|{\rm Im}\,\lambda|^{1/2}\|f\|_{1,0}+h^{1/2}\|f\|_{1,0}.$$
Taking $h$ small enough we absorb the last term and arrive at the estimate
\begin{equation}\label{eq:6.6}
\|f\|_{1,0}\le C|\lambda|^{1/2}|{\rm Im}\,\lambda|^{1/2}\|f\|_{1,0}
\end{equation}
for $\lambda\in\Lambda_k\setminus{\cal L}_N$, $0\le h\le h_0(N)$, with a constant $C>0$ independent of $\lambda$
and $N$. If 
\begin{equation}\label{eq:6.7}
|\lambda|^{-N}\le|{\rm Im}\,\lambda|\le (2C)^{-2}|\lambda|^{-1},\quad{\rm Re}\,\lambda\ge h_0(N)^{-1},
\end{equation}
 we can absorb the term in the right-hand side of
(\ref{eq:6.6}) and conclude that $\|f\|_{1,0}=0$, which implies $f\equiv 0$. In other words there are no
transmission eigenvalues in the region (\ref{eq:6.7}). It is easy to see that this implies that there are no
transmission eigenvalues in the region (\ref{eq:1.8}), provided the constant $C_N$ is properly chosen.

In what follows we will assume the condition (\ref{eq:1.9}) fulfilled and we will show that in this case
the factor $|\lambda|^{1/2}$ in the right-hand side of
(\ref{eq:6.6}) can be removed.  Clearly, the condition (\ref{eq:1.9}) implies ${\cal H}_2\subset {\cal H}_1$
and ${\cal E}_1\subset {\cal E}_2$. Given a parameter $0<\varepsilon\ll 1$, independent of $\lambda$, we can choose
functions $\chi_\varepsilon^-, \chi_\varepsilon^0, \chi_\varepsilon^+\in C^\infty(T^*\Gamma)$ such that 
$\chi_\varepsilon^-+\chi_\varepsilon^0+\chi_\varepsilon^+\equiv 1$, supp$\,\chi_\varepsilon^-\subset{\cal H}_1$,
supp$\,\chi_\varepsilon^+\subset{\cal E}_1$, $\chi_\varepsilon^0=1$ in $\{|cr_0/n-1|\le\varepsilon\}$ and 
$\chi_\varepsilon^0=0$ in $T^*\Gamma\setminus\{|cr_0/n-1|\le 2\varepsilon\}$.
Clearly, supp$\,\chi_\varepsilon^+\subset{\cal E}_2$. 
Taking $\varepsilon$ small enough we can also arrange that supp$(1-\chi_\varepsilon^-)\subset{\cal E}_2$.
Using this we will prove the following

\begin{lemma} We have the estimate
\begin{equation}\label{eq:6.8}
\|{\rm Op}_h(1-\chi_\varepsilon^-)f\|_{1,0}\lesssim (|{\rm Im}\,\lambda|^{1/2}+h^{1/2}+\varepsilon)\|f\|_{1,0}
\end{equation}
for $\lambda\in\Lambda_k\setminus{\cal L}_N$.
\end{lemma}

{\it Proof.} Observe that the condition (\ref{eq:1.9}) implies the inequality
$$\frac{c_2n_2-c_1n_1}{c_2^2-c_1^2}<\frac{n_1}{c_1}.$$
Therefore the inequality (\ref{eq:6.4}) holds on ${\cal E}_1$. This implies the estimate (\ref{eq:6.5}) with $1-\chi$ replaced by $\chi_\varepsilon^+$.
On the other hand, since the function $\chi_\varepsilon^+$ is supported in both elliptic regions,
the estimate (\ref{eq:6.3}) holds with $1-\chi$ replaced by $\chi_\varepsilon^+$. Thus we conclude that 
the estimate (\ref{eq:6.2}) holds with $1-\chi$ replaced by $\chi_\varepsilon^+$, that is, we have
\begin{equation}\label{eq:6.9}
\|{\rm Op}_h(\chi_\varepsilon^+)f\|_{1,0}\lesssim (|{\rm Im}\,\lambda|^{1/2}+h^{1/2})\|f\|_{1,0}.
\end{equation}
Choose a function  $\eta_\varepsilon\in C_0^\infty(T^*\Gamma)$
such that $\eta_\varepsilon=1$ in $\{|cr_0/n-1|\le 3\varepsilon\}$ and 
$\eta_\varepsilon=0$ in $T^*\Gamma\setminus\{|cr_0/n-1|\le 4\varepsilon\}$. Clearly, $\eta_\varepsilon=1$
on supp$\,\chi_\varepsilon^0$. Moreover, taking $\varepsilon$ small enough we can arrange that 
supp$\eta_\varepsilon\subset{\cal E}_2$. It is easy to see that $Tf=0$ implies the identity
$${\rm Op}_h(\chi_\varepsilon^0\rho_2)f+\left({\rm Op}_h(\chi_\varepsilon^0){\rm Op}_h(\eta_\varepsilon\rho_2)f-{\rm Op}_h(\chi_\varepsilon^0\rho_2)f\right)$$
$$={\rm Op}_h(\chi_\varepsilon^0){\rm Op}_h(\eta_\varepsilon\rho_2)f=-{\rm Op}_h(\chi_\varepsilon^0){\cal N}_2(\lambda){\rm Op}_h(\eta_\varepsilon)f$$
$$+{\rm Op}_h(\chi_\varepsilon^0)\left({\cal N}_2(\lambda){\rm Op}_h(\eta_\varepsilon)f+{\rm Op}_h(\eta_\varepsilon\rho_2)f\right)$$
$$={\rm Op}_h(\chi_\varepsilon^0){\cal N}_2(\lambda){\rm Op}_h(1-\eta_\varepsilon)f+{\rm Op}_h(\chi_\varepsilon^0)c_2^{-1}c_1{\cal N}_1(\lambda)f$$
$$+{\rm Op}_h(\chi_\varepsilon^0)\left({\cal N}_2(\lambda){\rm Op}_h(\eta_\varepsilon)f+{\rm Op}_h(\eta_\varepsilon\rho_2)f\right).$$
The $h-\Psi$DO calculus give
$${\rm Op}_h(\chi_\varepsilon^0){\rm Op}_h(\eta_\varepsilon\rho_2)-{\rm Op}_h(\chi_\varepsilon^0\rho_2)
={\cal O}_\varepsilon(h):L^2(\Gamma)\to L^2(\Gamma),$$
$${\rm Op}_h(\chi_\varepsilon^0)c_2^{-1}c_1-c_2^{-1}c_1{\rm Op}_h(\chi_\varepsilon^0)
={\cal O}_\varepsilon(h):L^2(\Gamma)\to L^2(\Gamma).$$
Therefore, using the estimates (\ref{eq:2.33}), (\ref{eq:4.1}), (\ref{eq:5.2}), (\ref{eq:5.14}) and (\ref{eq:6.1}), we get
\begin{equation}\label{eq:6.10}
\|{\rm Op}_h(\chi_\varepsilon^0\rho_2)f\|_0\lesssim (|{\rm Im}\,\lambda|^{1/2}+h^{1/2}+\varepsilon)\|f\|_{1,0}.
\end{equation}
On the other hand, the condition (\ref{eq:1.9}) guarantees that $|\rho_2|\ge C>0$ on supp$\,\chi_\varepsilon^0$,
provided $\varepsilon$ is taken small enough. Hence
\begin{equation}\label{eq:6.11}
\|{\rm Op}_h(\chi_\varepsilon^0)f\|_0\lesssim 
\|{\rm Op}_h(\chi_\varepsilon^0\rho_2)f\|_0+h\|f\|_0.
\end{equation}
Clearly, the estimate (\ref{eq:6.8}) follows from (\ref{eq:6.9}), (\ref{eq:6.10}) and (\ref{eq:6.11}).
\eproof

By (\ref{eq:2.27}), (\ref{eq:6.1}) and (\ref{eq:6.8}) we obtain 
$$\|f\|_{1,0}\lesssim |{\rm Im}\,\lambda|^{1/2}\|f\|_{1,0}+h^{1/2}\|f\|_{1,0}+\varepsilon\|f\|_{1,0}.$$
Taking $h$ and $\varepsilon$ small enough we can absorb the last two terms to obtain
\begin{equation}\label{eq:6.12}
\|f\|_{1,0}\le C|{\rm Im}\,\lambda|^{1/2}\|f\|_{1,0}
\end{equation}
for $\lambda\in\Lambda_k\setminus{\cal L}_N$, $0\le h\le h_0(N)$, with a constant $C>0$ independent of $\lambda$
and $N$. If 
\begin{equation}\label{eq:6.13}
|\lambda|^{-N}\le|{\rm Im}\,\lambda|\le (2C)^{-2},\quad{\rm Re}\,\lambda\ge h_0(N)^{-1},
\end{equation}
 we can absorb the term in the right-hand side of
(\ref{eq:6.12}) and conclude that $\|f\|_{1,0}=0$, which implies $f\equiv 0$. Hence there are no
transmission eigenvalues in the region (\ref{eq:6.13}), which implies that there are no
transmission eigenvalues in the region (\ref{eq:1.10}), provided the constant $C_N$ is properly chosen.


\begin{thebibliography}
\frenchspacing \baselineskip=12 pt plus 1pt minus 1pt

\bibitem{kn:BLR} {\sc C. Bardos, G. Lebeau and J. Rauch}, {\em Sharp sufficient conditions for the observation, control and stabilization
of waves from the boundary}, SIAM J. Control Optim. {\bf 30} (1992), 1024-1065.

\bibitem{kn:CCH} {\sc F. Cakoni, D. Colton and H. Haddar}, {\em Transmission eigenvalues}, Notices of the AMS {\bf 68} (9) (2021), 1499-1510.

\bibitem{kn:CPV} {\sc F. Cardoso, G. Popov and G. Vodev}, {\em Distribution of resonances and local energy decay in the transmission
problem. II}, Math. Res. Lett. {\bf 6} (1999), 377-396.

\bibitem{kn:G} {\sc J. Galkowski}, {\em The quantum Sabine law for resonances in transmission problems}, 
Pure Appl. Analysis {\bf 1} (1) (2019), 27-100.

\bibitem{kn:FN} {\sc J. Fornerod and H-M. Nguyen}, {\em The Weyl law of transmission eigenvalues and the completeness of generalized 
transmission eigenfunctions without complementing conditions}, SIAM J.
Math. Anal. {\bf 55} (2023), no. 4, 3959-3999.

\bibitem{kn:LV1} {\sc E. Lakshtanov and B. Vainberg}, {\em Ellipticity in the interior transmission problem in anisotropic media},
SIAM J. Math. Anal. {\bf 44} (2012), no. 2, 1165-1174.

\bibitem{kn:LV2} {\sc E. Lakshtanov and B. Vainberg}, {\em Applications of elliptic operator theory to the 
isotropic interior transmission eigenvalue problem}, Inverse Problems {\bf 29} (2013), 104003.

\bibitem{kn:LR} {\sc G. Lebeau and L. Robbiano}, {\em Contr\^ole exact de l'\'equation de la chaleur}, Commun. Partial Diff. Equations
{\bf 20} (1995), 335-356.

\bibitem{kn:NN1} {\sc H-M. Nguyen and Q-H. Nguyen}, {\em Discreteness of interior transmission eigenvalues revisited}, 
Calc. Var. Partial Differential Equations {\bf 56} (2017), no. 2, Paper No. 51.

\bibitem{kn:NN2} {\sc H-M. Nguyen and Q-H. Nguyen}, {\em The Weyl law of transmission eigenvalues and the completeness of generalized 
transmission eigenfunctions}, J. Funct. Anal. {\bf 281} (2021), 109146.

\bibitem{kn:PV} {\sc V. Petkov and G. Vodev}, {\em Asymptotics of the number of the interior transmission eigenvalues}, J. Spectral Theory
{\bf 7} (2017), 1-31. 

\bibitem{kn:R1} {\sc L. Robbiano}, {\em Spectral analysis of interior transmission eigenvalues},  Inverse Problems {\bf 29} (2013), 104001.

\bibitem{kn:R2} {\sc L. Robbiano}, {\em Counting function for interior transmission eigenvalues}, Math. Control Relat. Fields
{\bf 6} (2016), no. 1, 167-183.

\bibitem{kn:Sj} {\sc J. Sj\"ostrand}, {\em Singularit\'es analytiques microlocales}, Ast\'erisque {\bf 95} (1982). 

\bibitem{kn:Sy} {\sc J. Sylvester}, {\em Discreteness of transmission eigenvalues via upper triangular compact operators}, SIAM J.
Math. Anal. {\bf 44} (2012), no. 1, 341-354.

\bibitem{kn:V1} {\sc G. Vodev}, {\em Transmission eigenvalue-free regions}, Comm. Math. Phys. {\bf 336} (2015), 1141-1166.

\bibitem{kn:V2} {\sc G. Vodev}, {\em High-frequency approximation of the interior Dirichlet-to-Neumann map and applications to the
transmission eigenvalues}, Anal. PDE {\bf 11} (2018), 213-236.

\bibitem{kn:V3} {\sc G. Vodev}, {\em Interior transmission problems with coefficients of low regularity}, preprint 2023.




\end{thebibliography}
\end{document}